\documentclass[12pt,a4paper,reqno]{amsart}

\usepackage{amsmath}
\usepackage{amsfonts}
\usepackage{amssymb}
\usepackage{graphicx}

\usepackage[latin1]{inputenc}
\usepackage{amsmath}
\usepackage{amsfonts}
\usepackage{amsthm}
\numberwithin{equation}{section}
\newtheorem{teo}{Theorem}[section]
\newtheorem{lem}[teo]{Lemma}

\newtheorem{cor}[teo]{Corollary}
\theoremstyle{definition}
\newtheorem{defi}{Definition}
\theoremstyle{remark}
\newtheorem{rem}{Remark}

\newcommand{\ve}{\varepsilon}

\newcommand{\deb}{\rightharpoonup}

\newcommand{\pical}{\mathcal{P}}
\newcommand{\lcal}{\mathcal{L}}

\newcommand{\R}{\mathbb{R}}
\newcommand{\F}{\mathfrak{F}}
\newcommand{\tto}{\rightarrow}

\newcommand\huno{{\mathcal H}^1}

\newcommand{\N}{\mathbb{N}}

\DeclareMathOperator{\bari}{bar}
\DeclareMathOperator{\spt}{spt}

\author{F. Santambrogio}
\address{F. Santambrogio, Scuola Normale Superiore, Piazza dei Cavalieri 7, 56126 Pisa, Italy, {\tt f.santambrogio@sns.it}}
\title{Transport and Concentration Problems with Interaction Effects}

\begin{document}

\begin{abstract}
After presenting an overview about variational problems on probability measures for functionals involving transport costs and extra terms encouraging or discouraging concentration, we look for optimality conditions, regularity properties and explicit computations in the case where Wasserstein distances and interaction energies are considered.
\end{abstract}
\maketitle

\section{Introduction}\label{sec1}

One of the aim of the paper is to give a short overview of possible variational problems involving transport costs between distributions of mass and their concentrations. The general problem we are interested in is
$$\min_{\mu,\nu\in\pical(\Omega)} \F(\mu,\nu):= T(\mu,\nu)+F(\mu)+G(\nu),$$
where the functional $T$ quantifies in some way the distance between the two probability measures $\mu$ and $\nu$ according to a transport cost criterium, and $F$ and $G$ are functionals over the space $\pical(\Omega)$ (the space of probability measures over a domain $\Omega$) with opposite behaviour: the first prefers spread measures and penalizes concentration while the latter, on the other hand, prefers concentrated measures. Obviously there are lots of sub-problems which may be of interest, for instance the minimizations of the two separate functionals
$$\F_{\nu}(\mu):=T(\mu,\nu)+F(\mu)\quad\text{and}\quad\F^{\mu}(\nu)=T(\mu,\nu)+G(\nu),$$
where each time one of the variable is frozen. Also imposing constraints like $F(\mu)\leq H,\,G(\nu)\leq L\dots$ instead of adding penalizations in the functionals may sometimes be considered (and this is in fact the same as adding penalizations through some $0/+\infty$ functions).

These minimization problems are likely to appear in several phenomena both in nature and in decision science. For instance in \cite{ButSan}, \cite{CarSan} and \cite{mathese} these variational problems have been intended for urban planning models, where the spread measure $\mu$ stands for inhabitants, the concentrated $\nu$ for services and they have to be close in a transportation distance sense. Recently a possible choice of the functional $\F_{\nu}$ has been proposed as a model for the formation of a certain kind of leaves: if $\nu=\delta_0$ representes the point where feeding arrives to the leaf, then the shape of the leaf is such that it optimizes the quantity of light it receives from the sun taking also into account the effective trasport cost for bringing feeding all over its shape.

We present here some choices for the functionals $T$ and $G$. The choice of $F$ is in fact easier since a very good class of concentration-penalizing functionals is given by local convex functionals over measures, for instance
\begin{equation}\label{defiF}
F(\mu)=\begin{cases}\int_{\Omega}f(u)\,d\lcal&\text{ if }\mu=u\cdot\lcal\\
+\infty&\text{ otherwise,}\end{cases}
\end{equation}
for any convex fucntion $f$ with $f(0)=0$ which is superlinear at infinity.
For these functionals we refer to \cite{bb1}. Here $\lcal$ is a reference measure that may be chosen as the Lebesgue measure $\lcal^d$ if $\Omega\subset\R^d$. By Jensen inequality, spread measures with constant density are optimal for these functionals.

Possible choices of $T$ are the following:
\begin{itemize}
\item terms involving Monge-Kantorovich optimal transport cost, as in Wasserstein distances $T(\mu,\nu)=W_p^p(\mu,\nu)$ ({\it Wasserstein});
\item terms taking into account traffic congestion effects in trasportation distances, as in \cite{CarSan} where the choice $T(\mu,\nu)=||\mu-\nu||_{X'}^2$ for a vector space $X\subset H^1(\Omega)$ is detailed and justified ({\it congestion});
\item terms reflecting the natural ramified structure of a transportation network as in \cite{xia1}, where a new distance on probability measures is introduced according to this criterium ({\it branching}).
\end{itemize}
This last possibility is the most suitable for model involving natural branching structures like leaves, while the first two seem to be quite natural in urban planning. 

For the functional $G$, before presenting a list of examples, we give a possible definition of the concept of concentration-preferring:
\begin{defi}
We say that $G:\pical(\Omega)\tto\R$ is a {\it concentration-preferring} functional if it holds $G(t_{\sharp}\nu)\leq G(\nu)$ for any measure $\nu\in\pical(\Omega)$ and any $1-$Lipschitz continuous map $t:\Omega\tto\Omega$.
\end{defi}
It is easy to show that any $G$ with this property is minimized by any measure $\delta_{x_0}$, with $x_0\in\Omega$.
We list here some funcionals satisfying this definition:
\begin{itemize}
\item $G(\nu)=\sharp\spt(\nu)$ ({\it atomic}), as in location problems, where the corresponding $T$ is usually of Wasserstein type;
\item ({\it subadditive}, see \cite{bb1})
$$G(\nu) = \begin{cases}
                 \sum_{k\in\N}g(a_k) & \textrm{if} \ \nu=\sum_{k\in\N}a_k\delta_{x_k} \\
                 +\infty & \textrm{otherwise}    \end{cases} $$
            for a subadditive function $g$ with $g(0)=0$ and $g'(0)=+\infty$;
\item $G(\nu)=\inf\left\{\huno(\Sigma)\left|\spt(\nu)\subset\Sigma,\,\Sigma\text{ compact and connected }\right.\right\}$ as in the irrigation problem, \cite{ButOudSte}, where $T=W_1$ and a constraint on $G$ is considered ({\it length});
\item $G(\nu)=\int_{\Omega\times\Omega}w(d(x,y))\nu(dx)\nu(dy)$ for an increasing function $w$ ({\it interaction}).
\end{itemize}
Actually the two first terms are decreasing under the effect of any map $t$ and not only under $1-$Lipschitz ones.
The last choice for $G$ is in fact a well-known functional on probability measures which representes interaction energy (or costs) between the particles composing $\nu$. It has been first studied in a transportation framework by McCann in\cite{mccann}, where displacement convexity results are given, with the aim of producing uniqueness for variational problems.

In \cite{ButSan} and \cite{CarSan} two combinations of these functional have been studied for urban planning purposes: the {\it Wasserstein} + {\it subadditive} and the {\it congestion} + {\it interaction} cases, respectively. The {\it congestion} + {\it subadditive} case has been exluded in \cite{CarSan} since it leads to a trivially infinite functional, and so in this paper we analyze the remaining one, i.e. the {\it Wasserstein} + {\it interaction} case. Many ideas are taken from \cite{CarSan}, up to the fact that elliptic regularity is replaced by considerations on Monge-Ampère equation. Moreover some extra devices are performed and a careful use of Monge-Kantorovich theory is needed.

\section{Preliminaries on Optimal Transportation}\label{sec2}

We recall in this section the tools and definitions we need in the sequel about optimal transport theory and Monge-Ampère equation. Our main reference is \cite{villani}.

\begin{defi}
Given two probability measures $\mu$ and $\nu$ on a space $\Omega$ and a l.s.c. cost function $c:\Omega\times\Omega\tto[0,+\infty]$ we consider the problem
\begin{equation}\label{kantorovich}
(K)\quad\min\left\{\int_{\Omega\times\Omega}c\,d\gamma\left|\gamma\in\pical(\Omega\times\Omega),\,(\pi_1)_{\sharp}\gamma=\mu,\,(\pi_2)_{\sharp}\gamma=\nu\right.\right\},
\end{equation}
and the minimizers for this problem are called {\it optimal transport plan} between $\mu$ and $\nu$. Should $\gamma$ be of the form $(id\times t)_{\sharp}\mu$ for a measurable map $t:\Omega\tto\Omega$, the map $t$ would be called {\it optimal transport map} from $\mu$ to $\nu$.
\end{defi}

An important tool will be duality theory and to introduce it we need in particular the notion of $c-$transform (a kind of generalization af the well-known Legendre transform).

\begin{defi} 
Given a function $\chi:\Omega\tto\overline{\R}$ we define its {\it $c-$transform} (or $c-$conjugate function) by
$$\chi^c(y)=\inf_{x\in\Omega}c(x,y)-\chi(x).$$
Moreover, we say that a function $\psi$ is {\it $c-$concave} if there exists $\chi$ with $\psi=\chi^c$ and we denote by $\Psi_c(\Omega)$ the set of $c-$concave functions.
\end{defi}

It is well-known a duality result stating the following equality:
\begin{equation}\label{duality formula}
\min(K)=\sup_{\psi\in\Psi_c(\Omega)}\int_{\Omega}\psi\,d\mu+\int_{\Omega}\psi^c\,d\nu.
\end{equation}
\begin{defi} 
The functions $\psi$ realizing the maximum in \eqref{duality formula} are called {\it Kantorovich potentials} for the transport from $\mu$ to $\nu$.
\end{defi}
Since we will use $c(x,y)=\frac{1}{2}|x-y|^2$, let us denote by $\Psi_2(\Omega)$ the set of $c-$concave functions with respect to this quadratic cost. It is not difficult to check that 
$$\psi\in\Psi_2(\Omega) \Leftrightarrow \;x\mapsto \frac{x^2}{2}-\psi\text{  is a convex function on $\R^d$ restricted to $\Omega$.}$$
Notice that on a bouded $\Omega$ with diameter $D$ any $\psi\in\Psi_2(\Omega)$ is in fact $2D-$Lipschitz continuous.
We summarize here some useful results for the quadratic case $c(x,y)=\frac{1}{2}|x-y|^2$.
\begin{teo}
Given $\mu$ and $\nu$ probability measures on a connected $\Omega\subset\R^d$ there exists unique an optimal transport plan $\gamma$ and it is of the form $(id\times t)_{\sharp}\mu$, provided $\mu$ is absolutely continuous. Moreover there exists also at least a Kantorovich potential $\psi$, and the gradient $\nabla\psi$ is uniquely determined $\mu-$a.e. (in particular $\psi$ is unique up to additive constants, provided the density of $\mu$ is positive a.e. on $\Omega$). The optimal map $t$ and the potential $\psi$ are linked by $t(x)=x-\nabla\psi(x)$ and so $t$ is the gradient of a convex function. Moreover it holds $\psi(x)+\psi^c(t(x))=c(x,t(x))$ for $\mu-$a.e. $x$.
\end{teo}
Starting from the values of the problem $(K)$ in \eqref{kantorovich} we can define a set of distances over $\pical(\Omega)$. For any $p\geq 1$ we can define
$$W_p(\mu,\nu)=\big(\min(K)\text{ with }c(x,y)=|x-y|^p\big)^{1/p},$$
and obviously we will restrict our analysis to the case $p=2$. We recall that it holds, by duality formula,
\begin{equation}\label{duality formula2}
\frac{1}{2}W_2^2(\mu,\nu)= \sup_{\psi\in\Psi_2(\Omega)} \int_{\Omega}\psi\,d\nu + \int_{\Omega}\psi^c\,d\mu.
\end{equation}
\begin{teo}
If $\Omega$ is compact, for any $p\geq 1$ the function $W_p$ is in fact a distance over $\pical(\Omega)$ and the convergence with respect to this distance is equivalent to the weak convergence of probability measures. In particular any functional $\mu\mapsto W_p(\mu,\nu)$ is continuous with respect to weak topology.
\end{teo}

The next step of our analysis is concerned with some regularity properties of $t$ and $\psi$ (the optimal transport map and the Kantorovich potential, respectively) and their relations with the densities of $\mu$ and $\nu$. It is easy, just by a change-of-variables formula, to transfrom, in the case of regular maps and densities, the equality $\nu=t_{\sharp}\mu$ into the PDE $v(t(x))=u(x)/|Jt|(x)$, where $u$ and $v$ are the densities of $\mu$ and $\nu$ and $J$ denotes the determinant of the Jacobian matrix. Recalling $t=\nabla\phi$ with $\phi$ convex, we get the Monge-Ampère equation
\begin{equation}\label{mongeampere}
M\phi=\frac{u}{v(\nabla\phi)},
\end{equation}
where $M$ denotes the determinant of the Hessian $M\phi=\det H\phi= \det \left[\frac{\partial^2\phi}{\partial x_i\,\partial x_j}\right]_{i,j} $.
This equation up to now is satisfied by $\phi=id-\psi$ just in a formal way. We define various notion of solutions for \eqref{mongeampere}:
\begin{itemize}
\item we say that $\phi$ satisfies \eqref{mongeampere} in the Brenier sense if $(\nabla\phi)_{\sharp}u\cdot\lcal^d=v\cdot\lcal^d$ (and this is actually the sense to be given to this equation);
\item we say that $\phi$ satisfies \eqref{mongeampere} in the Aleksandrov sense if $H\phi$, which is always a positive measure for $\phi$ convex, is absolutely continuous and its density satisfies \eqref{mongeampere} a.e.;
\item we say that $\phi$ satisfies \eqref{mongeampere} in the viscosity sense if it satisfies the usual comparison properties required by viscosity theory but restricting the comparisons to regular convex test functions (since $M$ is in fact monotone just when restricted to positively definite matrices);
\item we say that $\phi$ satisfies \eqref{mongeampere} in the classical sense if is $C^2$ and the equation holds pointwise.
\end{itemize}
Notice that any notion but the first may be also applied to the equation $M\phi=f$, while the first one just applies to this specific transportation case.
The results we want to use are well summarized in Theorem 50 in \cite{villani}:
\begin{teo}
If $u$ and $v$ are $C^{0,\alpha}(\Omega)$ and are both bounded above and below on the whole $\Omega$ by positive constants and $\Omega$ is a convex open set, then for the unique Brenier solution $\phi$ of \eqref{mongeampere} it holds $\phi\in C^{2,\alpha}(\Omega)\cap C^{1,\alpha}(\overline{\Omega})$ and $\phi$ satisfies the equation in the classical sense (hence also in the Aleksandroff and viscosity senses).
\end{teo}



\section{Optimality Conditions for the Interaction Case}\label{sec3}

We are here concerned with the minimization problem for the functional $\F^{\mu}$, when the transport term is given by $T(\mu,\nu)=\frac{1}{2}W_2^2(\mu,\nu)$ and the concentration one is an interaction term of the form
\begin{equation}\label{defiG}
G(\nu)=\int_{\Omega\times\Omega}V(|x-y|^2)\nu(dx)\nu(dy),
\end{equation}
with $V:[0,+\infty[\tto[0,+\infty[$ a regular increasing fuction. From now on $\Omega$ will be the closure of a convex non-empty open set in $\R^d$. 

A priori a minimizer for this functional may be an arbitrary probability measure on the set $\Omega$, even a singular one. Our goal is to prove, under suitable assumptions and by means of optimality conditions and of an approximation process, that it is in fact an absolutely continuous measure with bounded density.

We provide here an easy optimality condition for the minimization of $\F^{\mu}$. We do not go into details in the computation both beacuse it is not difficult to perform it and because it follows the same, trivial, scheme as in \cite{ButSan}. Notice that the approximation argument by measures with positive densities that we are going to use works in this case too, while the alternative proof by convex analysis that may be found in \cite{ButSan} does not, simply because there is no convexity in the term $G$.
\begin{teo}\label{optimality}
If a probability measure $\nu\in\pical(\Omega)$ is a minimizer for $\F^{\mu}$, then there exists a constant $m$ such that
$$\psi + T_{\nu} \geq m; \quad \psi + T_{\nu} = m \quad\nu\text{-a.e.} ,$$
where $\psi$ is a Kantorovich potential for the transport from $\nu$ to $\mu$ and we define
$$T_{\nu}(x)=\int_{\Omega}2V(|x-y|^2)\,\nu(dy).$$
\end{teo}
\begin{proof}
Let us start from the case when $\mu$ is absolutely continuous with positive density. In this case we perturbate an optimal measure $\nu$ into $\nu_t=\nu + t (\nu_1-\nu)$ for an arbitrary $\nu_1\in\pical(\Omega)$: if we call $\psi_t$ the unique Kantorovich potential from $\nu_t$ to $\mu$ which vanishes at a certain fixed point $x_0\in\Omega$, we get (by means of duality formula)
\begin{multline*}
\int_{\Omega}\psi_t\,d\nu_t+\int_{\Omega}\psi_t^c\,d\mu+\int_{\Omega\times\Omega}V(|x-y|^2)\nu_t(dx)\nu_t(dy)\\
\geq
\int_{\Omega}\psi_t\,d\nu+\int_{\Omega}\psi_t^c\,d\mu+\int_{\Omega\times\Omega}V(|x-y|^2)\nu(dx)\nu(dy).
\end{multline*}
After erasing the term $\int_{\Omega}\psi_t^c\,d\mu$ and dividing by $t$ we pass to the limit, and we know that $\psi_t$ converges uniformly (by Ascoli-Arzelà) to the unique Kantorovich potential $\psi$ from $\nu$ to $\mu$ vanishing at $x_0$. This provides, at the limit,
$$\int_{\Omega}(\psi+T_{\nu})\,d\nu_1\geq\int_{\Omega}(\psi+T_{\nu})\,d\nu.$$
Being $\nu_1$ arbitrary we get that $\nu-$a.e. the function $\psi+T_{\nu}$ must be equal to its infimum, and this is the thesis.

To generalize the result to an arbitrary measure $\mu$, just proceed by approximation. This can be performed standardly as in \cite{ButSan} and provides the same formula where $\psi$ becomes just one of the possibly many Kantorovich potentials instead of the only one.
\end{proof}

The problem in this condition lies in the fact that the measure $\nu$ appears only in a very implicit way (both in $\psi$ and in $T_{\nu}$), and this does not allow to derive any estimate on it. We will consequently need to pass through an aproximation process, exactly as in \cite{CarSan}. Fixed a minimizer $\bar{\nu}$ for $\F^{\mu}$, we will consider a sequence of problems $(P_{\ve})_{\ve}$ given by the minimization of
$$\pical(\Omega)\ni\nu\mapsto T(\mu_{\ve},\nu)+ G(\nu)+\delta_{\ve}A(\nu)+\ve W_2^2(\nu,\bar{\nu}_{\ve}),$$
where
\begin{itemize}
\item $(\mu_{\ve})_{\ve}$ is a sequence of probability measures approximating $\mu$ with Lipschitz continuous strictly positive densities $u_{\ve}$;
\item the functional $A$ is given by 
$$A(\nu)=\begin{cases}\int_{\Omega}a(v)\,d\lcal^d&\text{ if }\nu=v\cdot\lcal^d,\\
                      +\infty&\text{ otherwise,}\end{cases}$$
for a convex function $a:[0,+\infty[\tto[0,+\infty]$ which is superlinear at infinity and blowing up at $0$, i.e. $\lim_{t\tto 0^+}a(t)=+\infty$, but finite and $C^2$ with $a''\geq c>0$ on $]0,+\infty[$ (for instance $a(t)=t^2+1/t$);
\item $(\delta_{\ve})_{\ve}$ is a suitably chosen sequence of small positive numbers;
\item $(\bar{\nu}_{\ve})_{\ve}$ is a suitably chosen sequence of measures approximating $\bar{\nu}$.
\end{itemize}
We will prove, exactly as in \cite{CarSan}, that this sequence of problems admits an uniform $L^{\infty}$ bound for their solutions and that we can choose the parameter so that these solutions converge to $\bar{\nu}$, thus obtaining an $L^{\infty}$ estimate for $\bar{\nu}$. The existence of solutions for $P_{\ve}$ is trivial by the semicontinuity of each term in the sum with respect to the weak convergence of probability measures on the compact set $\Omega$ (and moreover any term but $A$ is actually continuous, while $A$ is semicontinuous by convexity).

\begin{lem}\label{uniform bound}
Suppose that $\mu$ is absolutely continuous with an $L^{\infty}$ density bounded by $M$ and that $V$ is a $C^2$ function with $V'>0$, then any solution $\nu_{\ve}$ to the problem $P_{\ve}$ is absolutely continuous and its density is bounded by a constant $C$ depending only on $M,\,d$ and $V$.
\end{lem}
\begin{proof}
First we notice that, thanks to the presence of the term $A(\nu)$ in the minimization problem, $\nu_{\ve}$ must be absolutely continuous with strictly positive density almost everywhere. We write down the optimality conditions for $\nu_{\ve}$ with respect to variations of the form $\nu_t=\nu_{\ve}+t(\nu_1-\nu_{\ve})$. From easy computations we get
$$\psi_{\ve}+T_{\nu_{\ve}}+\delta_{\ve}a'(\nu_{\ve})+\ve\chi_{\ve}=m_{\ve}\text{  q.o.},$$
where $\psi_{\ve}$ is the Kantorovich potential for the tansport from $\nu_{\ve}$ to $\mu_{\ve}$ and $\chi_{\ve}$ from $\nu_{\ve}$ to $\bar{\nu}_{\ve}$ (they are unique up to additive constants) and $m_{\ve}$ is a suitable constant. We get equality almost everywhere due to the fact that we already know that $\nu_{\ve}>0$ (we identify measures and their densities in this context). Since Kantorovich potentials are Lipschitz functions and $T_{\nu_{\ve}}$ shares the same regulartiy of the integrand $(x,y)\mapsto V(|x-y|^2)$, which is $C^2$ and then Lipschitz, we get that even $a'(\nu_{\ve})$ is Lipschitz continuous, and in particular it is bounded. This prevents $\nu_{\ve}$ to be close to $0$ since it holds $\lim_{t\tto 0^+}a'(t)=-\infty$. Thus we get $\nu_{\ve}\geq c_{\ve}>0$. Moreover, $a'(\nu_{\ve})$ is Lipschitz continuous and, being $a''$ strictly greater than a positive constant, also the inverse of $a'$ is Lipschitz. This proves that $\nu_{\ve}$ is a Lipschitz continuous function.
We can now use regularity theory on Monge-Ampère equation to get $\psi\in C^{2,\alpha}(\Omega)\cap C^{1,\alpha}(\overline{\Omega})$, since both $\nu_{\ve}$ and $\mu_{\ve}$ are bounded both from above and from below by positive finite constants (depending on $\ve$, anyway) and are Lipschitz continuous.
The same is true for the Kantorovich potential $\chi_{\ve}$ by replacing $\mu_{\ve}$ by $\bar{\nu}_{\ve}$. 
What we can do now is looking for a maximum point $x_0$ of $\nu_{\ve}$. Notice that such a point will be a minimum point for $\psi_{\ve}+T_{\nu_{\ve}}+\ve\chi_{\ve}$. First we prove that $x_0\notin\partial\Omega$. To prove this it is sufficient to prove that the gradient of $\psi_{\ve}+T_{\nu_{\ve}}+\ve\chi_{\ve}$ is directed outwards at any point of $\partial\Omega$, i.e. $\nabla(\psi_{\ve}+T_{\nu_{\ve}}+\ve\chi_{\ve})(x_0)\cdot n(x_0)>0$ for any $x_0\in\partial\Omega$, where $n$ is the outward normal vector. From the fact that the optimal transport map $t$  from $\nu_{\ve}$ to $\mu_{\ve}$ is given by $t(x)=x-\nabla\psi(x)$ we know that $x-\nabla\psi(x)\in\Omega$ for almost any $x\in\Omega$. In this case, due to continuity up to the boundary of $\nabla\psi$, this holds for any $x$ and also for $x_0\in\partial\Omega$ and implies $\nabla\psi(x_0)\cdot n(x_0)\geq 0$. Analogously we get $\nabla\chi(x_0)\cdot n(x_0)\geq 0$. For the gradient of $T_{\nu_{\ve}}$ it holds
$$\nabla T_{\nu_{\ve}}(x_0)=\int_{\Omega}4V'(|x_0-y|^2)(x_0-y)\,\nu_{\ve}(dy),$$
and so $\nabla T_{\nu_{\ve}}(x_0)\cdot n(x_0)>0$ since $V'>0$ and for almost any $y\in\Omega$ it holds $(x_0-y)\cdot n(x_0)> 0$.
This proves that $x_0$ lies in the interior of $\Omega$ and this allows us to look at the second derivatives. Taking Hessians we have
$$H\psi_{\ve}(x_0)+HT_{\nu_{\ve}}(x_0)+\ve H\chi_{\ve}(x_0)\geq 0,$$
where the letter $H$ denotes Hessians and the inequality is in the sense of positive definite symmetric matrices.
Thus we get 
$$H\psi_{\ve}(x_0)\geq -I \left(2||V||_{C^2(\Omega)}+\ve \right),$$
since the second derivatives of $T_{\nu_{\ve}}$ may be estimated by those of $V$ and from the fact that $x^2/2-\chi(x)$ is convex we deduce $H\chi\leq I$. This is a uniform estimate from below for $H\psi_{\ve}(x_0)$ and for the convex function $\phi$ given by $\phi(x)=x^2/2-\psi_{\ve}(x)$ we have $H\phi(x_0)\leq I \left( 1+\ve + 2||V||_{C^2(\Omega)}\right).$ This implies $M\phi(x_0)\leq ( 1+\ve + 2||V||_{C^2(\Omega)})^d$, and, from $\nu_{\ve}=\mu_{\ve}(\nabla \phi) M\phi$, we get, for $\ve\leq 1$,
$$\max \,\nu_{\ve}=\nu_{\ve}(x_0)\leq 2^dM \left( 1 + ||V||_{C^2(\Omega)}\right)^d,$$
which is the desired estimate.
\end{proof}

\begin{rem}
The simple proof we give here of the fact that the gradient is directed outwards and no maximum point is allowed on the boundary could be used similarly in \cite{CarSan}, thus getting rid of the strict convexity assumption in Theorem 6.5 and of the heavy proof of Lemma 6.6. Notice that it could be possible to get the same result even without $C^1$ regularity for the potentials, just making the proof a bit trickier. It would be sufficient to evaluate the increments of the potential in small balls around $x_0$ where the gradient is almost everywhere defined and such that $x-\nabla\psi(x),\;x-\nabla\chi(x)\in\Omega$ for a.e. $x$.
\end{rem}

\begin{lem}\label{convergence is possible}
It is possible to choose the parameters for the problem $P_{\ve}$, i.e. the numbers $\delta_{\ve}$ and the measures $\bar{\nu}_{\ve}$ and $\mu_{\ve}$ so that any sequence of minimizers $(\nu_{\ve})_{\ve}$ converges to $\bar{\nu}$.
\end{lem}
\begin{proof}
It is sufficient to choose for $\bar{\nu}_{\ve}$ a sequence of absolutely continuous measures with Lipschitz continuous strictly positive densities such that $F(\bar{\nu}_{\ve})\leq F(\bar{\nu})+\ve^2$. Then we have $A(\bar{\nu}_{\ve})<+\infty$ and we may choose $\delta_{\ve}=\ve^2 A(\bar{\nu}_{\ve})^{-1}$. For $(\mu_{\ve})_{\ve}$ we can chose any sequence of absolutely continuous measures with Lipschitz continuous strictly positive densities approximating $\mu$ in such a way that $W_2(\mu_{\ve},\mu)\leq \ve^2$. Then we have
\begin{eqnarray*}
T(\mu_{\ve},\nu_{\ve})\!+\!G(\nu_{\ve})\!+\!\delta_{\ve}A(\nu_{\ve})\!+\!\ve W_2^2(\nu_{\ve}, \bar{\nu}_{\ve})\leq T(\mu_{\ve},\bar{\nu}_{\ve})\!+\!G(\bar{\nu}_{\ve})\!+\!\delta_{\ve}A(\bar{\nu}_{\ve}),
\end{eqnarray*}
which implies
\begin{eqnarray*}
F(\nu_{\ve})+\delta_{\ve}A(\nu_{\ve})+\ve W_2^2(\nu_{\ve}, \bar{\nu}_{\ve})&\leq& F(\bar{\nu}_{\ve})+4DW_2(\mu_{\ve},\mu)+\delta_{\ve}A(\bar{\nu}_{\ve})\\
&\leq&  F(\bar{\nu}) + \ve^2 +4D\ve^2 + \ve^2\\
&\leq& F (\nu_{\ve})+\ve^2(2+4D).
\end{eqnarray*}
In the end this implies $W_2(\nu_{\ve},\bar{\nu}_{\ve})\leq C\sqrt{\ve}$ and, since $\bar{\nu}_{\ve}\deb\bar\nu$, we get $\nu_{\ve}\deb\bar{\nu}$.
\end{proof}

We can now state our main result and its consequence in the minimization for the whole functional $\F$.

\begin{teo}\label{princ}
Given a compact convex set $\Omega\subset\R^d$ with nonempty interior and a probability measure $\mu\in L^{\infty}(\Omega)$, if the function $V$ appearing in the definition of the functional $G$ is $C^2$ and $V'>0$, then the minimization problem for the functional $\F^{\mu}$ over the space $\pical(\Omega)$ admits at least a solution and any minimizer belongs in fact to the space $L^{\infty}(\Omega)$.
\end{teo}
\begin{proof}
As usual the existence is trivial due to continuity and compactness of $\pical(\Omega)$ while, for the $L^{\infty}$ regularity, just apply Lemma \ref{uniform bound} and Lemma \ref{convergence is possible}
\end{proof}

\begin{cor}
Given a compact convex set $\Omega\subset\R^d$ with nonempty interior, a $C^1$ strictly convex and superlinear function $f$ and a $C^2$ function $V$ with $V'>0$, then the minimization problem over the space $\pical(\Omega)^2$ for the functional $\F(\mu,\nu)= \frac{1}{2}W_2^2(\mu,\nu)+F(\mu)+G(\nu)$, where $F$ is defined by \eqref{defiF} and $G$ by \eqref{defiG} admits solutions and, in any minimizing pair $(\mu,\nu)$, both $\mu$ and $\nu$ are in fact absolutely continuous measures $\mu=u\cdot\lcal^d,\,\nu=v\cdot\lcal^d$, with $u\in C^0(\Omega)$ and $v\in L^{\infty}(\Omega)$.
\end{cor}
\begin{proof}
After the usual proof of existence by the direct method in Calculus of Variations, we refer to \cite{ButSan} for the regularity results on $\mu$. Since such a measure turns out to be absolutely continuous with continuous density (hence bounded), we may apply Theorem \ref{princ} to get the regularity on $\nu$.
\end{proof}

\section{An Explicit Example}\label{sec4}

In this section we come back to the whole problem of minimizing $\F$ in a very particular case, where we can provide almost explicit densities fo the solutions. We consider the case
\begin{itemize}
\item $T(\mu,\nu)=\frac{1}{2} W_2^2(\mu,\nu)$ and $G(\nu)=\int_{\Omega\times\Omega}V(|x-y|^2)\nu(dx)\nu(dy),$ as in the previous Section;
\item $V(|x-y|^2)=\frac{\lambda}{2}|x-y|^2$ and so, setting $\bari(\nu)=\int_{\Omega}y\nu(dy)$, we have $T_{\nu}(x)=\lambda|x|^2-2\lambda x\cdot\bari(\nu) + \lambda \int_{\Omega}|y|^2\nu(dy)$;
\item $F(\mu)=\frac{1}{2}||\mu||_{L^2(\Omega)}^2,$ a particular case of what considered in \cite{ButSan}.
\end{itemize}
The framework we obtain is very similar to the one in \cite{CarSan}.

\begin{teo}
In the specific case described above, any pair of minimizers $(\mu,\nu)$ is shaped as follows:
\begin{itemize}
\item $\mu$ is concentrated on a ball $B(x_0,r)$ (intersected with $\Omega$) and has a density $u$ given by $u(x)=\frac{\lambda}{2\lambda+1}(r^2-|x-x_0|^2)$;
\item $\nu$ is concentrated on the ball $B(x_0,r/(2\lambda+1))$ and 
it is the image of $\mu$ under the omothety of center $x_0$ and ratio $(2\lambda+1)^{-1}$;
\item $x_0$ is the barycenter of $\nu$.
\end{itemize}
\end{teo}
\begin{proof}
First we write down the optimalty conditions given by Theorem \ref{optimality} for the minimization in $\nu$ with fixed $\mu$ and by \cite{ButSan} for the minimization in $\mu$ for fixed $\nu$. We denote by $u$ and $v$ the densities of $\mu$ and $\nu$, respectively. We may suppose that the barycenter of $\nu$ is the origin, thus obtaining $T_{\nu}(x)=\lambda|x|^2+c$. We have
$$\begin{cases} u(x)+\varphi(x) =c_1&\text{ a.e. on }u>0;\\
                \psi(x)+\lambda x^2= c_2&\text{ a.e. on }v>0.\end{cases}$$
Here $\varphi$ and $\psi$ are Kantorovich potentials for the transport from $\mu$ to $\nu$ and from $\nu$ to $\mu$, respectively.
From the second condition we can infer $\nabla\psi(x)=-2\lambda x$ a.e. on $v>0$. Being $\nu$ absolutely continuous, this equality is valid $\nu-$a.e.. This means that the optimal transport map $t$ from $\nu$ to $\mu$ is given by $t(x)=x-\nabla\psi(x)=(2\lambda+1)x$. By uniqueness of the optimal transport plan, which is in this case expressed both as a transport map from $\nu$ to $\mu$ and viceversa, we know also the optimal transport map $s$ from $\mu$ to $\nu$ which will be $s(x)=x/(2\lambda+1)$.
From duality theory in optimal transportation we know the following equality
$$\varphi(x)+\varphi^c(s(x))=c(x,s(x))=\frac{1}{2}|x-s(x)|^2,$$
and thus we get $u(x)=c_1-\frac{1}{2}|x-s(x)|^2+\varphi^c(s(x))$. Since $\varphi^c$ is a Kantorovich potential from $\nu$ to $\mu$, we know that it must agree (up to constants) with $\psi$ on any connected component of $\{v>0\}$ (notice that in this case we are allowed to speak about connected components of this set, since $\{u>0\}$ is an open set, being $u$ continuous, and the set $\{v>0\}$ is just an omothety of it). So, let $\omega\subset\Omega$ be a connected component of $\{u>0\}$. On $(2\lambda+1)^{-1}\omega$ we have $\varphi^c=\psi+c_3$ and so we get
$$u(x)=c_4-\frac{1}{2}|x-s(x)|^2+\psi(s(x))=
c_5-|x|^2\frac{\lambda}{2\lambda+1}.$$
From this expression it is clear that $\partial\omega\setminus\partial\Omega$ (where $u$ must vanish) is contained in a circle line around $0$. This implies that $0$ belongs in fact to $\omega$, since no boundary of $\omega$ is allowed in the interior of a certain ball around $0$. So there is in fact just one connected component for $\{u>0\}$ and so we get
\begin{equation}\label{expliu}
u(x)=\left[c-|x|^2\frac{\lambda}{2\lambda+1}\right]^+.
\end{equation}
From this it is easy to recover the density of $\nu$ since $\nu=s_{\sharp}\mu$ and we get the thesis. The point $x_0$ wich turns out to be the center of the balls which are supports for $\mu$ and $\nu$ is in this notation $0$, the barycenter of $\nu$, as in the thesis. It is clear that in this case $\mu$ and $\nu$ share the same barycenter since they are omothetical. 
\end{proof}

\begin{rem}
In this example the density $v$ shares the same regularity of $u$ except at the points corresponding to boundary points of $\Omega$ where $u$ is positive, i.e. if at $x_0\in\partial\Omega$ it happens $u(x_0)>0$ then at $s(x_0)$ we have a jump for $v$. It is clear from the fact that $u$ is $2\lambda/(2\lambda+1)-$Lipschitz continuous (it follows from the explicit formula) that we have, recalling also $\int_{\Omega}u\,d\lcal^d=1$,
$$1\leq\left(\inf u+\frac{2\lambda}{2\lambda+1}D\right)|\Omega|,$$
where $D$ is the diameter of $\Omega$. This implies, for small $\Omega$, $\inf u>0$. In this case $u$ would be positive at any point of $\partial\Omega$ and $v$ discontinuous at any point of $s(\partial\Omega)$. This gives examples when the $L^{\infty}$ regularity for $v$ is quite sharp (in the sense $v\notin C^0(\Omega)$).
\end{rem}
\begin{rem}
In the explicit example above there remain to be determined both the constant $c$ in the formula for $u$ and the position of the barycenter $x_0$. In some simple cases this too is possible. Notice that, once fixed $x_0$, the constant $c$ may always be recovered by imposing the condition of being probability measures. For instance if $\Omega$ is a ball, we may see that $x_0$ may not be the barycenter of a density $u$ shaped as in \eqref{expliu} unless the set $B(x_0,2\lambda^{-1}(2\lambda+1))\cap\Omega$ is a ball around $x_0$. This happens for large $\Omega$ whenever the ball $B(x_0,2\lambda^{-1}(2\lambda+1))$ does not touch the boundary $\partial\Omega$ or, in general, when $x_0$ is the center of the ball $\Omega$. In the first case ($\Omega$ a large ball) we have several solutions for the problem (non-uniqueness), obtained from each other under translations, and $u$ and $v$ are continuous; in the second ($\Omega$ a small ball) we have uniqueness of the solution, with $u$ a radial continuous function around the center and $v$ a rescaling of $u$ on a smaller ball.
\end{rem}


\end{document}